\title{ ~~\\ Improvement of
an estimate of H. M\"uller involving the order of $2({\rm mod~}u)$ II}
\author{Pieter Moree}
\documentclass[12pt]{article}
\usepackage{amssymb, latexsym, amsfonts}
\def\@ptsize{2}
\newtheorem{Thm}{Theorem}

\newtheorem{Lem}{Lemma}
\newtheorem{Cor}{Corollary}

\newtheorem{Prop}{Proposition}
\newcommand{\qed}{\hfill $\Box$}

\begin{document}
\date{}
\maketitle
\begin{abstract}
\noindent Let $m\ge 1$ be an arbitrary fixed integer and 
let $N_m(x)$ count the number of odd integers $u\le x$ such
that the order of $2$ modulo $u$ is not divisible by $m$.
In case $m$ is prime estimates for $N_m(x)$ were given by M\"uller 
that were subsequently sharpened into 
an asymptotic estimate by the present author. M\"uller
on his turn extended the author's result to the case where $m$ is
a prime power and gave bounds in the case $m$ is not a prime
power. Here an asymptotic for $N_m(x)$ is derived that is valid for
all integers $m$. This asymptotic would easily have followed from 
M\"uller's approach were it
not for the fact that a certain Diophantine equation has non-trivial
solutions. All solutions of this equation are determined. We also
generalize to other base numbers than $2$. For a very sparse set
of these numbers M\"uller's approach does work.
\end{abstract}
\section{Introduction}
Let $u$ be odd. Denote by $l(u)$ the smallest natural number such
that $$2^{l(u)}\equiv 1({\rm mod~}u).$$ The number $l(u)$ is called
the {\it order} of the congruence class $2({\rm mod~}u)$.
Let $N_m(x)$ denote the number of odd integers $u\le x$ such that 
$m\nmid l(u)$. It was shown by Franco and Pomerance \cite[Theorem 5]{FP} in 
the context of a generalization of the celebrated "$3x+1$"-problem, that almost 
all integers
$u$ have the property that $m|l(u)$, i.e. they established that, 
as $x$ tends to infinity,
$$N_m(x)=o(x).$$
The object of this note is to derive an asymptotic formula for
$N_m(x)$. Partial progress towards this goal was made by H. 
M\"uller in his papers \cite{HM1, HM2}. 
Let $q>2$ be a prime. M\"uller \cite{HM1} showed
that 
$${x\over \log^{1/(q-1)}x}\ll N_q(x)\ll {x\over \log^{1/q}x}.$$
This was improved by the present author in \cite{PM1}, where
he showed that
\begin{equation}
\label{pietjebell}
N_q(x)=c_q{x\over \log^{q/(q^2-1)}x}\left(1+O\left({(\log \log x)^5\over 
\log x}\right)\right),
\end{equation}
with $c_q>0$ a positive real constant.
On his turn M\"uller \cite{HM2} improved on this, by establishing the
following generalization of the asymptotic estimate (\ref{pietjebell}). 
\begin{Thm} {\rm (H. M\"uller)}.
\label{mullertwee}
Let $q>2$ be a prime and $n\ge 1$ fixed. Then
$$N_{q^n}(x)=c_{q^n}{x\over \log^{q^{2-n}/(q^2-1)}x}\left(1+O\left({(\log \log x)^5\over 
\log x}\right)\right),$$
with $c_{q^n}$ a positive real constant.
\end{Thm}
M\"uller did not obtain an asymptotic
result in case $m$ is not a prime power. If 
$m=p_1^{e_1}\cdots p_r^{e_r}$, he notes that
\begin{equation}
\label{smullie}
N_{p_1^{e_1}}(x)\le N_m(x)\le N_{p_1^{e_1}}(x)+\cdots+N_{p_r^{e_r}}(x)
\end{equation}
and from this infers, in case $m$ is odd, that
$${x\over \log^{\alpha} x}\ll N_m(x)\ll {x\over \log^{\beta} x},$$
where
$\alpha:={\rm max}\{C(p_j,e_j)|1\le j\le r\},~{\rm resp.~}
\beta:={\rm min}\{C(p_j,e_j)|1\le j\le r\}$, 
where for any prime $q$ and integer 
$n\ge 1$ we define $C(q,n):=q^{2-n}/(q^2-1)$.
(In fact, M\"uller erroneously swaps `min' and `max' in his definitions of
$\alpha$ and $\beta$.)\\
\indent Note that the stronger result, in case $m$ is odd, that
$x\log^{-\beta}x\ll N_m(x)\ll x\log^{-\beta}x$ actually follows from
(\ref{smullie}). More can be said however:
\begin{Lem}
\label{blok}
Let $m\ge 1$ be odd. Suppose that there is only one integer $j$ such that
$C(p_j,e_j)=\beta$, then we have that
$N_m(x)\sim N_{p_j^{e_j}}(x)$ as $x$ tends to infinity and the 
asymptotic behaviour of
$N_m(x)$ is given by Theorem {\rm \ref{mullertwee}}.
\end{Lem}
{\it Proof}. W.l.o.g. we can assume that $C(p_1,e_1)=\beta$. Now if there is
no further $j$ such that that $C(p_j,e_j)=\beta$, then by Theorem \ref{mullertwee} each of the
terms $N_{p_k^{e_k}}(x)$ with 
$k\ge 2$ is asymptotically of smaller growth than
$N_{p_1^{e_1}}(x)$  and so the lemma 
is proved. \qed\\

Note that if there is more than one integer $j$ such 
that $C(p_j,e_j)=\beta$,
then from (\ref{smullie}) together with 
Theorem \ref{mullertwee} the asymptotic behaviour of $N_m(x)$ cannot
be determined. 
We are thus led to the problem of determining pairs of integers $j$ and
$k$ and primes $p_j$ and $p_k$ such that $C(p_j,e_j)=C(p_k,e_k)$. 
The proof of the following result uses some arguments kindly provided
by J.-H. Evertse and Y. Bilu. It follows that $C(2,5)=C(5,2)=C(3,3)=1/24$, $C(2,6)=C(7,2)=1/48$ 
and $C(5,3)=C(11,2)=1/120$.
\begin{Prop}
There are only finitely many solutions 
$(p,\alpha,q,\beta)$ to 
\begin{equation}
\label{gelijk}
p^{\alpha-2}(p^2-1)=q^{\beta-2}(q^2-1)
\end{equation}
with
$\alpha$ and $\beta$ integers and $p$ and $q$ primes 
with $p<q$. These solutions are: 
$(2,5,3,3)$, $(2,5,5,2)$, $(2,6,7,2)$, $(3,3,5,2)$ and
$(5,3,11,2)$.
\end{Prop}
{\it Proof}. It is easy to see that if there is to be a solution, we
must have $\alpha\ge 2$ and $\beta\ge 2$. 
Assume w.l.o.g. that $q>p$. If $\beta\ge 3$, then $q$
must divide either $p-1$ or $p+1$. Since $q>p$ it follows that
$q=p+1$ and so $p=2$ and $q=3$. This gives rise to precisely one
solution: $(2,5,3,3)$. So we may assume that $\beta=2$ and we
are reduced to finding the solutions of 
\begin{equation}
\label{gelijk2}
p^{b}(p^2-1)=q^2-1,
\end{equation}
with $b\ge 0$.\\
\indent First assume that $p=2$. Then we have to solve the
equation $q^2-1=3\cdot 2^b$. As one of $q\pm 1$ cannot be divisible by
$4$, it follows that either $q-1\in \{1,2,3,6\}$ or 
$q+1\in \{1,2,3,6\}$. This gives rise
to the solutions $(2,5,5,2)$ and $(2,6,7,2)$ (and no more). 
Thus we may assume that $p>2$. There are two cases to be considered:\hfil\break
1) $p^b|q+1$. Then we can write $q+1=p^b r$, $q-1=s$ (say) and
so $rs=p^2-1$. Thus $p^br-s=2$ and hence $r(p^br-2)=p^2-1$. Note that 
if $b>1$, then we must have $r>1$ and this gives rise to a contradiction.
So $b=1$ and hence we must have $r(pr-2)=p^2-1$. From this we infer
that $0<r<p$ and, furthermore that $-2r\equiv -1({\rm mod~}p)$. It
follows that $r=(p+1)/2$. Substituting this back we obtain 
$p^2-2p-3=0$ which only gives rise to the solution $(3,3,5,2)$.\\
2) $p^b|q-1$. Now we have that $q+1$ divides $p^2-1$ and we obtain
$p^b\le q-1<q+1\le p^2-1$, which is impossible when $b>1$.
We can write $q-1=pr$, $q+1=s$ (say) and
so $rs=p^2-1$. Thus $s-pr=2$ and hence $r(pr+2)=p^2-1$.   
From this we infer
that $0<r<p$ and, furthermore that $-2r\equiv 1({\rm mod~}p)$. It follows
that $r=(p-1)/2$. Substituting this back we see
that $p^2-6p+5=0$ which only gives rise to the solution $(5,3,11,2)$.\\ 
\indent On collecting all solutions found along the way, the proof is 
completed. \qed\\

\noindent Since the Diophantine equation (\ref{gelijk}) has non-trivial
solutions we are blocked for certain integers $m$ (for example for
$m=2^5\cdot 5^2=800$), in proving an asymptotic for $N_m(x)$ by
using (\ref{smullie}). We will provide an 
alternative to (\ref{smullie}), Lemma \ref{burp}, 
and use this to obtain a more precise estimate for $N_m(x)$ than
can be provided by M\"uller's method, and works for every integer $m$.
At the same time we generalize to other base numbers than $2$.

\section{Generalization to other base numbers}
First we will generalize the main result mentioned sofar, Theorem 1, to
the case where the base number $g$ is rational and not in $\{-1,0,1\}$ (an 
assumption on $g$ maintained throughout this paper). 
We let $\omega(n)$ denote the number of distinct prime divisors of $n$
and by $\nu_p(n)$ denote the exponent of $p$ in the prime factorization 
of $n$. Let $S(g)$ be the set of integers composed only of primes $p$ that
do not occur in the prime factorization of $g$ (i.e. of primes $p$ such that $\nu_p(g)=0$). For each 
integer
$u$ in this set the order of $g$ modulo $u$, ord$_g(u)$, is well-defined. 
We let $P_g(d)(x)$ denote the number of primes $p\le x$ with $p\in S(g)$ such
that $d|{\rm ord}_g(p)$. It turns out that the set $P_g(d)$ of primes $p\in S(g)$ such
that $d|{\rm ord}_g(p)$ has a natural density, which 
will be denoted by $\delta_g(d)$. 
This density was first
determined by Wiertelak \cite{W}, who 
derived a rather complicated explicit
formula for it, a formula which was subsequently simplified
by Pappalardi \cite{Pappa} and streamlined further by Moree \cite{Polen}. It turns
out that $\delta_g(d)$ is always a positive rational number.
Part 1 of the
following result is due to Wiertelak \cite{W}, part 2 to Moree \cite{Polen}. As usual
the logarithmic integral is denoted by Li$(x)$.
\begin{Thm} \label{wimo} {\rm (Wiertelak \cite{W}, Moree \cite{Polen})}.\\
\noindent {\rm 1)} We have $P_g(d)(x)=\delta_g(d){\rm Li}(x)+O\left({x\over \log^3 x}(\log \log x)^{\omega(d)+3}\right)$.\\
{\rm 2)} Under GRH we have $P_g(d)(x)=\delta_g(d){\rm Li}(x)+O(\sqrt{x}\log^{\omega(d)+1} x )$.
\end{Thm}
\noindent In order to explicitly evaluate $\delta_g(d)$, some further
notation is needed. Write $g=\pm g_0^h$, where $g_0$ is
positive and not an
exact power of a rational and $h$ as large as possible. 
Let $D(g_0)$ denote the discriminant of the field $\mathbb Q(\sqrt{g_0})$.
(This notation will reappear several times in the sequel.)
Given an integer $d$, we denote by $d^{\infty}$ the supernatural number (sometimes
called Steinitz number), $\prod_{p|d}p^{\infty}$. Note 
that ${\rm gcd}(v,d^{\infty})=\prod_{p|d}p^{\nu_p(v)}$.\\

\noindent {\tt Definition}. Let $d$ be even and let $\epsilon_g(d)$ be defined as in Table 1 with
$\gamma={\rm max}\{0,\nu_2(D(g_0)/dh)\}$.

\centerline{{\bf Table 1:} $\epsilon_g(d)$}
\medskip
\begin{center}
\begin{tabular}{|c|c|c|c|}\hline
$g\backslash \gamma$ & $\gamma=0$&$\gamma=1$&$\gamma=2$\\
\hline\hline
$g>0$&$-{1/2}$&${1/4}$&${1/16}$\\
\hline
$g<0$&${1/4}$&$-{1/2}$&${1/16}$\\
\hline
\end{tabular}
\end{center}
\medskip

\noindent Note that $\gamma\le 2$. Also note that $\epsilon_g(d)=(-1/2)^{2^{\gamma}}$ if $g>0$.
\begin{Thm}
\label{main} {\rm (Explicit evaluation of $\delta_g(d)$)}.
We have
$$\delta_g(d)={\epsilon'_1(d)\over d{\rm ~gcd}(h,d^{\infty})}\prod_{p|d}{p^2\over p^2-1},{\rm ~with~}$$
$$\epsilon'_1(d)=\cases{1 & if $2\nmid d$;\cr
1+3(1-{\rm sgn}(g))(2^{\nu_2(h)}-1)/4& if $2||d$ and $D(g_0)\nmid 4d$;\cr
1+3(1-{\rm sgn}(g))(2^{\nu_2(h)}-1)/4+\epsilon_g(d)& if $2||d$ and $D(g_0)|4d$;\cr
1 & if $4|d$, $D(g_0)\nmid 4d$;\cr
1+\epsilon_{|g|}(d) & if $4|d$, $D(g_0)|4d$.}
$$
\end{Thm}
\indent By $N(x;g,m)$ we denote the number of 
integers $u\le x$ such that $u\in S(g)$ and
$m\nmid {\rm ord}_g(u)$. Note that $N(x;2,m)=N_m(x)$.  
A straightforward generalization of Theorem 1 yields the following result.
\begin{Thm} 
\label{mullerdire}
Let  $n\ge 1$ be fixed. Then
$$N(x;g,q^n)=c_{q^n}(g){x\over \log^{\delta_g(q^n)}x}\left(1+O\left({(\log \log x)^5\over 
\log x}\right)\right).$$
\end{Thm}
It is an easy consequence of Theorem \ref{main} that
$\delta_g(q^n)=C(q,n)$ for almost all integers $g$. Furthermore, by Theorem \ref{main}
we find that
$$\delta_2(q^n)=\cases{17/24 & if $q=2$ and $n=1$;\cr
5/12 & if $q=2$ and $n=2$;\cr
2^{1-n}/3 & if $q=2$ and $n\ge 3$;\cr
C(q,n) & otherwise.}
$$
This evaluation allows us to formulate Theorem 1 for every prime power (and thus
also for every power of two). It then follows from (2) that 
$x\log^{-\gamma_2(m)}x\ll N_m(x)\ll x\log^{-\gamma_2(m)}x$, where
$$\gamma_g(m):=\min \{\delta_g(p_j^{e_j})~|~1\le j\le r\}.$$ 
Lemma 1
can then be extended to all natural numbers $m$ as follows:
\begin{Lem}
\label{blok2}
Suppose that there is only one integer $j$ such that
$\delta_2(p_j^{e_j})=\gamma_2(m)$, then we have that
$N_m(x)\sim N_{p_j^{e_j}}(x)$ as $x\rightarrow \infty$ and the 
asymptotic behaviour of
$N_m(x)$ follows from Theorem {\rm \ref{mullerdire}}.
\end{Lem}
We leave it as an exercise to the reader to show that the density of integers
$m$ satisfying such that $\delta_2(p_j^{e_j})=\gamma_2(m)$ for only one 
integer $j$ exists and equals
$${147497571941\over 147916692000}\approx 0.9971665\cdots$$
An integer $g$ such that there is only one integer $j$ such that
$\delta_g(p_j^{e_j})=\gamma_g(m)$ for every natural number $m$, we
define to be a {\it M\"uller number}. M\"uller's inequality
(2) generalizes of course to
\begin{equation}
\label{ongelijkzeg}
N(x;g,p_1^{e_1})\le N(x;g,m)\le N(x;g,p_1^{e_1})+\cdots + N(x;g,p_r^{e_r}).
\end{equation}
The usefulness of M\"uller numbers is apparent from the following result.
\begin{Prop}
If $g$ is a M\"uller number and $\delta_g(p_j^{e_j})=\gamma_g(m)$, then 
we have $N_m(x)\sim N_{p_j^{e_j}}(x)$ as $x\rightarrow \infty$ and the 
asymptotic behaviour of
$N_m(x)$ is given by {\rm Theorem \ref{mullerdire}}.
\end{Prop}
In the next section it will be seen that, unfortunately, M\"uller numbers
are very sparse.

\section{Some further Diophantine considerations}
The above discussion shows that if for $g=2$ we 
want to cover odd integers $m$ as well, rather than asking for the
solutions of $C(p_i,e_i)=C(p_j,e_j)$, we should be asking
for solutions of $\delta_2(p_i^{e_i})=\delta_2(p_j^{e_j})$, and
indeed, more generally, for the non-trivial solutions 
$(p_1,e_1,p_2,e_2)$ with $p_1\ne p_2$ primes and 
$e_1,e_2\ge 1$ of
\begin{equation}
\label{algemener}
\delta_g(p_1^{e_1})=\delta_g(p_2^{e_2}).
\end{equation}
A variation of the proof of Proposition 1 gives that in case $g=2$ we
only have the non-trivial equalities $\delta_2(2^4)=\delta_2(3^3)=\delta_2(5^2)$, 
$\delta_2(2^5)=\delta_2(7^2)$ and $\delta_2(5^3)=\delta_2(11^2)$.
The more general situation is described by the following result.
\begin{Thm}
\label{vier}
There are only finitely many 
quadruples $(p_1,e_1,p_2,e_2)$ with $p_1<p_2$ primes and $e_1,e_2\ge 1$ such
that $\delta_g(p_1^{e_1})=\delta_g(p_2^{e_2})$.
These quadruples are given as follows (and there are no further ones):\\
$$\cases{(2,5-\tau_1-\nu_2(h),3,3-\nu_3(h)) & if $2^{5-\tau_2}\nmid h$ and $3^3\nmid h$;\cr
(2,5-\tau_1-\nu_2(h),5,2-\nu_5(h)) & if $2^{5-\tau_2}\nmid h$ and $5^2\nmid h$;\cr
(2,6-\tau_1-\nu_2(h),7,2-\nu_7(h)) & if $2^{6-\tau_2}\nmid h$ and $7^2\nmid h$;\cr
(3,3-\nu_3(h),5,2-\nu_5(h)) & if $3^3\nmid h$ and $5^2\nmid h$;\cr
(5,3-\nu_3(h),11,2-\nu_{11}(h)) & if $5^3\nmid h$ and $11^2\nmid h$,}$$
where $$\tau_1(g)=\cases{1 & if $D(g_0)=8$;\cr
0 & otherwise,}{\rm ~and~}\tau_2(g)=\tau_1(g)+{1-{\rm sgn}(g)\over 2}.$$
The associated values of $\delta_g$ are $1/24,1/24,1/48,1/24$, respectively $1/120$.
\end{Thm}
\begin{Cor}
Let 
$$S_1=\{8000,165373,193600,196000,209088,4002075,4743200,5122656\}.$$
Let $S_2$, respectively $S_4$ be $S_1$, but with all even numbers in this
set divided by $2$, respectively $4$.
For a given number $g$ equation {\rm (\ref{algemener})} has only trivial solutions, that
is $g$ is a M\"uller number, if and
only if $h$ is divisible by a number from $S_{2^{\tau_2(g)}}$.
\end{Cor}
\noindent {\tt Example}. The numbers $-2^{2000},2^{4000}$ and $3^{8000}$ are examples of
small M\"uller numbers. The number $g=3^{4000}$ is not a M\"uller number, since
$\delta_g(2^{1})=\delta_g(7^2)$ by Theorem \ref{vier}.\\

\noindent For the convenience of the reader we include a table where
$\tau_1=\tau_1(g)$ and $\tau_2=\tau_2(g)$ are given.\\

\centerline{{\bf Table 2:} $\tau_1(g)$ and $\tau_2(g)$}
\medskip
\begin{center}
\begin{tabular}
{|c|c|c|c|c|}\hline
&\multicolumn{2}{| c |}{$D(g_0)\ne 8$} & \multicolumn{2}{c |}{$D(g_0)=8$}\\ \hline
&$\tau_1$&$\tau_2$&$\tau_1$&$\tau_2$\\ \hline
$g>0$&$0$&$0$&$1$&$1$\\
\hline
$g<0$&$0$&$1$&$1$&$2$\\
\hline
\end{tabular}
\end{center}
\medskip

Clearly M\"uller numbers 
are very sparse. Indeed it is not difficult to quantify this: there are positive constants $c_{+}$ and $c_{-}$ 
such that the number of positive M\"uller numbers up to $x$ grows as $c_{+}x^{1/8000}$ and the number
of negative  M\"uller numbers with absolute value not exceeding $x$ as $c_{-}x^{1/4000}$ as $x\rightarrow 
\infty$.\\
\indent The proof of Theorem \ref{vier} is unfortunately not very instructive as many
cases have to be distinguished. As its level of difficulty is comparable to that of
the proof of Proposition 1, we leave it to the interested reader.\\
\indent In order
to derive Corollary 1 we notice that the set of those integers $h$ such
that none of the 5 solutions in Theorem \ref{vier} occur, are the
natural numbers of the form $(h_1)\cup \cdots \cup (h_s)$, where $(h_i)$
denotes the ideal generated by $h_i$ and $h_1,\ldots,h_s$ are certain integers. A priori the 5 solutions yield
$2^5=32$ potential generators 
(for example ${\rm lcm}(3^3,5^2,7^2,5^2,11^2)=4002075$), of which some turn out to be divisors of others (the latter
ones can thus be left out from $S_{2^{\tau_2(g)}}$).

\section{A more refined estimate for $N(x;g,m)$}
We sharpen M\"uller's estimate (\ref{ongelijkzeg}) to an equality (given in Lemma \ref{burp}) and use
this to obtain a more refined estimate for $N(x;g,m)$ (and thus $N_m(x)$).
As before we let $m={p_1}^{e_1}\cdots {p_r}^{e_r}$ denote the 
factorization of $m$. By $\kappa(m)=p_1\cdots p_r$ 
we denote the squarefree kernel of $m$. A divisor $d$ of $m$ is
said to be {\it unitary} if gcd$(d,m/d)=1$. Note that the 
non-trivial unitary divisors of $m$ come out as `blocks' from its 
factorization.
Let $N'(x;g,m)$ denote the number of integers $u\le x$ from $S(g)$ such
that $p_j^{e_j}\nmid {\rm ord}_g(u)$ for $1\le j\le r$. 
We let
$P'(x;g,m)$ be similarly defined, but with the phrase `integers $u\le x$'
replaced by `primes $p\le x$'. Note that $N'(x;g,p_j^{e_j})=N(x;g,p_j^{e_j})$.
\begin{Lem} 
\label{burp}
We have
$$N(x;g,m)=-\sum_{d\# m\atop d>1}\mu(\kappa(d))N'(x;g,d),$$
where by $d\# m$ we indicate that the sum is over the unitary divisors $d$
of $m$.
\end{Lem}
{\it Proof}. If $u\not\in S(g)$, the integer $u$ is not counted on either side of the purported identity. 
If $u\le x$ is in $S(g)$ we will show that it is counted
with multiplicity one on both sides of the purported identity. If $m|{\rm ord}_g(u)$, 
then $u$ is counted neither on the left nor on the right hand side.
If $m\nmid {\rm ord}_g(u)$, we can assume w.l.o.g. that ${p_j}^{e_j}\nmid {\rm ord}_g(u)$ for 
$1\le j\le k$ and that $p_s^{e_s}|{\rm ord}_g(u)$ for $s>k$. 
By definition $u$ is counted with multiplicity one on the left hand side.
The contribution of $u$
to the right hand side is 
$$-\sum_{d\# {p_1}^{e_1}\cdots {p_k}^{e_k}\atop d>1}\mu(\kappa(d))=
-\sum_{d|p_1\cdots p_k\atop d>1}\mu(d)=\mu(1)=1,$$
where we used the well-known identity $\sum_{d|n}\mu(d)=0$, which holds 
for every integer $n>1$. \qed\\

\noindent {\tt Example}. We have $N(x;2,12)=N_{12}(x)=N_3(x)+N_4(x)-N'(x;2,12)$.\\

Put $\delta'_g(d)=\sum_{j\# m}\mu(\kappa(j))\delta_g(j)$.
\begin{Lem} 
\label{lemma4}
We have\\
{\rm 1)} We have $P'(x;g,d)=\delta'_g(d){\rm Li}(x)+O\left({x\over \log^3 x}(\log \log x)^{\omega(d)+3}\right)$.\\
{\rm 2)} Under GRH we have $P'(x;g,d)=\delta'_g(d){\rm Li}(x)+O(\sqrt{x}\log^{\omega(d)+1} x )$.
\end{Lem}
{\it Proof}. By inclusion and exclusion (or an argument as in the proof of 
Lemma \ref{burp}) we see that
$P'(x;g,d)=\displaystyle\sum_{j\# d}\mu(\kappa(j))P_g(j)(x)$. Then use Theorem \ref{wimo}. \qed\\

\noindent {\tt Remark}. Wiertelak \cite{W} studied $P'(x;g,d)$ in depth in case $d$ is squarefree. 
Note (as did Wiertelak) that if $d$ is squarefree, then $P'(x;g,d)$ counts the number of primes $p\le x$ 
with $\nu_p(g)=0$ such that the congruence $g^{dy}\equiv g({\rm mod~}p)$ has a solution $y$. In
general  $P'(x;g,d)$ counts the number of primes $p\le x$ 
with $\nu_p(g)=0$ such that the congruence $g^{dy}\equiv g^r({\rm mod~}p)$ has a solution $(r,y)$ with $r$
dividing $d/\kappa(d)$.\\

The following result gives some insight in the properties of $\delta'_g(d)$.
\begin{Prop}
\label{pijnover}
Suppose that $d\ge 1$ and $d_1>1$ are natural numbers.\\
{\rm 1)} We have $\delta_g(dd_1)<\delta_g(d)$. If $g>0$, then we 
have $\delta_g(dd_1)\le {5\over 6}\delta_g(d)$.\\
{\rm 2)} We have 
$$\sup\{{\delta_g(dd_1)\over \delta_g(d)}|d_1>1\}=1{\rm ~and~}\sup\{{\delta_g(dd_1)\over \delta_g(d)}|g>0,~d_1>1\}=
{5\over 6}.$$
{\rm 3)} We have $\delta'_g(p_1^{e_1}p_2^{e_2})<\delta'_g(p_1^{e_1})$, with $p_1$ and $p_2$ primes.\\
{\rm 4)} We have $\delta'_g(p_1^{e_1}\cdots p_r^{e_r})=(1-\delta_g(p_1^{e_1}))\cdots (1-\delta_g(p_r^{e_r}))$, with
$p_1,\ldots,p_r$ distinct odd prime numbers.\\
{\rm 5)} We have $\min\{1-\delta'_g(j)|j\# d,~d>1\}=\gamma_g(d)$.
\end{Prop}
{\it Proof}. 1) If $g>0$, then one checks that $\epsilon'_1(dd_1)\le {5\over 4}\epsilon'_1(d)$ and 
from this one easily infers that $\delta_g(dd_1)\le {5\over 6}\delta_g(d)$. The general case requires 
some case distinctions and is left to the interested reader.\\
2) We have, e.g., $\lim_{e\rightarrow \infty}\delta_{-5^{2^e}}(6)/\delta_{-5^{2^e}}(3)=1$ and
$\delta_3(6)=5\delta_3(3)/6$. Now invoke part 1.\\
3) On invoking part 1 we find that
$\delta'_g(p_1^{e_1}p_2^{e_2})-\delta'_g(p_1^{e_1})=\delta_g(p_1^{e_1}p_2^{e_2})-\delta_g(p_2^{e_2})<0$.\\
4) Considered as a function of $d$, $\delta_g(d)$ is a multiplicative function on the odd integers
by Theorem \ref{main}. From this and the expression of $\delta'_g$ in terms of $\delta_g$, the
result then follows.\\
5) On noting that $\delta'_g(d d_1)\le \delta'_g(d)$ and using part 3 we see that
$$\min\{1-\delta'_g(j)|j\#d,~j>1\}=\min\{1-\delta'_g(p_i^{e_i})|1\le j\le r\}=\gamma_g(d),$$
where we used that $1-\delta'_g(p_i^{e_i})=\delta_g(p_i^{e_i})$. \qed\\

Next we will compute $N'(x;g,d)$ from $P'(x;g,d)$. To this end we will make use
of the fact that the integers counted by $N'(x;g,d)$ as $x\rightarrow \infty$ can be 
related to completely multiplicative sets.
Recall that a set $S$ of natural numbers is said to be 
{\it completely multiplicative} if its characteristic function $\chi$ is completely
multiplicative, that is satisfies $\chi(xy)=\chi(x)\chi(y)$ for all natural numbers $x$ and $y$.
Thus $S$ is completely multiplicative when  $x\cdot y\in S$
iff $x\in S$ and $y\in S$. The following result is a simple
generalization of Lemma 1 of M\"uller \cite{HM2} 
and easily follows from the observation that
${\rm ord}_g(u_1u_2)|{\rm gcd}(u_1,u_2){\rm lcm}({\rm ord}_g(u_1),{\rm ord}_g(u_2))$, 
cf. \cite{HM2}.

\begin{Lem}
\label{lemma5}
Let $m={p_1}^{e_1}\cdots {p_r}^{e_r}$. 
The set of integers $u$ from $S(g)$ such that $gcd(u,m)=1$ and $p_j^{e_j}\nmid {\rm ord}_g(u)$
for $1\le j\le r$ is a completely multiplicative set.
\end{Lem}

Let $N''(x;g,m)$ denote the number of integers $u\le x$ from $S(g)$ such
that gcd$(u,m)=1$ and $p_j^{e_j}\nmid {\rm ord}_g(u)$ for $1\le j\le r$. 
The integers thus counted (as $x\rightarrow \infty$) form a completely multiplicative set by Lemma \ref{lemma5} and can be
counted using the following proposition.
\begin{Prop} \label{pietjepuk}
{\rm 1)} Let $S$ be a completely multiplicative set of natural numbers such that 
for $x\rightarrow \infty$ we have
\begin{equation}
\label{tralala}
\sum_{p\le x\atop p\in S}1=\tau {\rm Li}(x)+O\left({x(\log \log x)^{\gamma}\over \log^3 x}\right)
\end{equation}
holds, where $\tau>0$ and $\gamma \ge 0$, with $\tau$ real and $\gamma$ an integer. Then we have
for $x\rightarrow \infty$ that
$$S(x):=\sum_{n\le x\atop n\in S}1=c_S x \log^{\tau-1}x + O(x(\log \log x)^{\gamma+1}\log^{\tau-2}x),$$
where $c_S$ denotes a positive constant which only depends on the set $S$.\\
{\rm 2)} If {\rm (\ref{tralala})} holds with error term $O(x\log^{-2-\gamma_1}x)$ for some $\gamma_1>0$, 
then $$S(x)=x\sum_{0\le \nu<\gamma_1}b_{\nu}\log^{\tau-1-\nu} x + O(x \log^{\tau-1-\gamma_1+\epsilon}x),$$
where $b_0(=c_S),b_1,\ldots $ are constants.
\end{Prop} 
{\it Proof}. The first part is proved in Moree \cite{Poud}. The second part is a special case of
Theorem 6 of \cite{MC}.\qed\\

Let  $m\ge 1$ be fixed. Then, by the latter proposition and Lemmas \ref{lemma4} and \ref{vier},
\begin{equation}
\label{doorduwen}
N''(x;g,m)=c'_{m}(g){x\over \log^{1-\delta'_g(m)}x}\left(1+O\left({(\log \log x)^{\omega(m)+4}\over 
\log x}\right)\right),
\end{equation}
where $c'_{m}(g)$ is a positive constant.\\
\indent As in \cite{HM2} we conclude that there exists a finite set of numbers (actually prime
powers) $n_1,\ldots,n_s$ (depending on $g$ and $m$) such that
\begin{equation}
\label{peelof}
N'(x;g,m)=\sum_{j=1}^s N''({x\over n_j};g,m).
\end{equation}
It thus follows on invoking the estimates (\ref{doorduwen}) and (\ref{peelof}) that
\begin{equation}
\label{doubleprime} 
N'(x;g,m)=c_{m}(g){x\over \log^{1-\delta'_g(m)}x}\left(1+O\left({(\log \log x)^{\omega(m)+4}\over 
\log x}\right)\right).
\end{equation}
On noting that $N'(x;g,p_j^{e_j})=N(x;g,p_j^{e_j})$ and that
$\delta'_g(p_j^{e_j})=1-\delta_g(p_j^{e_j})$, we see that the estimate 
(\ref{doubleprime}) generalizes Theorem \ref{mullerdire}.

Using (\ref{doubleprime}) we can now generalize Proposition 2.
\begin{Prop} 
Assume that $m=p_1^{e_1}\cdots p_r^{e_r}$ with
$\delta_g(p_1^{e_1})\le \cdots \le \delta_g(p_r^{e_r})$. Then we
have
$N(x;g,m)\sim \sum_{j=1}^{\min\{r,3\}}N(x;g,p_i^{e_i})$. It follows that in particular we have
$N(x;g,m)\sim c x\log^{-\delta_g(p_1^{e_1})}x$ for some $c>0$ as $x\rightarrow \infty$.
\end{Prop}
{\it Proof}. By Lemma \ref{burp}, (\ref{doubleprime}) and part 5 of Proposition \ref{pijnover}, we have
$N(x;g,m)\sim \sum_{j=1}^r N(x;g,p_i^{e_i})$. Since by Theorem \ref{vier} we
cannot have that $\delta_g(p_1^{e_1})=\delta_g(p_2^{e_2})=\delta_g(p_3^{e_3})=\delta_g(p_4^{e_4})$, 
the result follows on invoking (\ref{doubleprime}) again. \qed\\  

The main result of this paper makes the latter result more precise:
\begin{Thm}
Let $\{1-\delta'_g(j)~|~j\#m,~j>1\}=\{\delta_1,\ldots,\delta_s\}$ with $\delta_1<\delta_2<\cdots <\delta_s$.\\
{\rm 1)} Then $\delta_1=\gamma_g(m)$ and
$$N(x;g,m)=\sum_{j=1}^s c_j{x\over \log^{\delta_j}x}+O\left({x(\log \log x)^{5}\over \log^{1+\gamma_g(m)}x}
\right),$$
where the leading (asymptotically dominating) term is the one with $j=1$ and $c_1,\ldots,c_s$ are constants 
with $c_1>0$.\\
{\rm 2)} Under GRH it is true that for each integer $t\ge 0$ we have
$$N(x;g,m)=\sum_{j=1}^s \sum_{k=0}^t c_{j,k}{x\over \log^{k+\delta_j} x}+O({x\over \log^{t+1} x}).$$
The constants $c_j$ and $c_{j,k}$ will depend on $g$ and $m$ in general and so do the implied
constants in the error terms.
\end{Thm}
{\it Proof}. 1) Follows on combining Lemma \ref{burp}, (\ref{doubleprime}) and part 5 of 
Proposition \ref{pijnover}.\\
2) The proof is similar to that of part 1, except that instead of part 1 of Lemma \ref{lemma4} we use
part 2 and that instead of part 1 of Proposition \ref{pietjepuk} we use part 2.\qed\\

\noindent {\tt Acknowledgement}. The author thanks Yuri Bilu and
Jan-Hendrik Evertse for their assistance in solving the 
Diophantine equation (\ref{gelijk2}). The proof given is a compilation
of their arguments. Moreover, he thanks Igor Shparlinski for pointing
out the existence of \cite{HM2} to him.

\begin{small}
{\noindent Max-Planck-Institut f\"ur Mathematik\\ 
Vivatsgasse 7\\ 
D-53111 Bonn\\ 
Deutschland\\
E-mail:~moree@mpim-bonn.mpg.de\\}
\end{small}

\begin{thebibliography}{9}
\begin{small}
\bibitem{FP} Z. Franco and C. Pomerance, On a conjecture of Crandall
concerning the $qx+1$ problem, {\it Math. Comp.} {\bf 64} (1995), 1333-1336.
\bibitem{Poud} P. Moree, On the divisors of $a^k+b^k$, {\it Acta Arith.} {\bf 80} (1997), 197--212.
\bibitem{PM1} P. Moree, Improvement of an estimate of H. M\"uller 
involving the order of $2({\rm mod~}u)$, {\it Arch. Math.} {\bf 71} 
(1998), 197--200. 
\bibitem{Polen} P. Moree, On primes $p$ for which $d$ divides ord$_p(g)$, 
{\it Funct. Approx. Comment. Math.} {\bf 33} (2005), 85--95.
\bibitem{MC} P. Moree and J. Cazaran, On a claim of Ramanujan in his
first letter to Hardy, {\it Expo. Math.} {\bf 17} (1999), 289--312.
\bibitem{HM1} H. M\"uller, Eine Bemerkung \"uber die Ordnungen von
$2({\rm mod~}U)$ bei ungeradem $U$, {\it Arch. Math.} {\bf 69} (1997), 
217--220.
\bibitem{HM2} H. M\"uller, On the distribution of the orders
of $2({\rm mod~}u)$ for odd $u$, {\it Arch. Math.} {\bf 84} (2005), 
412--420.
\bibitem{Pappa} F. Pappalardi, Square free values of the order function, 
{\it New York J. Math.} {\bf 9}  (2003), 331--344. 
\bibitem{W} K. Wiertelak, On the density of some sets of primes. IV, 
{\it Acta Arith.} {\bf 43} (1984), 177--190.
\end{small}
\end{thebibliography}
\end{document}